\documentclass[12pt]{article}
\usepackage{amssymb,amsfonts,amsthm}
\usepackage{amstext}

\usepackage{amsmath}

\usepackage{color}
\usepackage[latin1]{inputenc}
\def\t{\theta}

\def\p{\phi}
\def\openC{{\rm C\kern-.18cm\vrule width.8pt height 7pt depth-.2pt \kern.18cm}}
\def\openN{{{\rm I}\kern-.16em {\rm N}}}
\def\openR{{{\rm I}\kern-.16em {\rm R}}}
\def\openT{{{\rm T}\kern-.42em {\rm T}}}
\def\openZ{{{\rm Z}\kern-.28em{\rm Z}}}

\def\eop{\hfill\rule{2.5mm}{2.5mm}}
\def\pf{\par\smallbreak\noindent {\bf Proof.} }

\newtheorem{thm}{Theorem}[section]

\newtheorem{lem}[thm]{Lemma}
\newtheorem{prop}[thm]{Proposition}

\theoremstyle{definition}

\def\eop{\hfill\rule{2.5mm}{2.5mm}}

\begin{document}

\title{\textbf{Strictly positive definite kernels on a product of circles} \vspace{-2pt}
\author{\sc
J. C. Guella\thanks{All authors partially supported by FAPESP, under grants $\#$ 2012/22161-3, $\#$2014/00277-5 and 2014/25796-5 respectively.}\, , V. A. Menegatto \,\,\&\,\,A. P. Peron}}
\date{}
\maketitle \vspace{-30pt}
\bigskip
\begin{center}
\parbox{13 cm}{{\small We supply a Fourier characterization for the real, continuous, isotropic and strictly positive definite kernels on a product of circles.}}
\end{center}
\bigskip
\noindent{\bf Key words and phrases:} strictly positive definite kernels, isotropy, product of circles, Schoenberg theorem, Skolem-Mahler-Lech theorem.\\
\noindent{\bf 2010 Math. Subj. Class.:}  33C50; 33C55; 42A16; 42A32; 42A82; 42B05; 43A35.

\thispagestyle{empty}

%
%
\section{Introduction}

Positive definite functions and kernels have a long history in mathematics, entering as an important tool in harmonic analysis and other areas as well.\ In the spherical setting, they can be traced back to
the remarkable paper of I. J. Schoenberg published in 1942 (\cite{schoen}), where a characterization for the continuous, isotropic and positive definite kernels on a single sphere was obtained.\ This characterization is far-reaching, having applications in approximation theory, spatial statistics, geomathematics, discrete geometry, etc.\ We mention \cite{cheney, dai,gneiting,musin} and references therein for some applications of positive definite functions and kernels on spheres.

In this paper, we will be concerned with positive definite kernels on a product of circles.\ As so, we will recast the basic concepts and results from Schoenberg's work that applies to circles, up to the point we can state what the main contribution in the present paper is.

We will write $S^1$ to denote the unit circle in $\mathbb{R}^{2}$.\ Continuity of a kernel $K$ on $S^1$ will be attached to the usual geodesic distance on $S^1$ and that will be extended to the product $S^1 \times S^1$ in the usual way.\ The {\em isotropy or radiality} of a kernel $K$ on $S^1$ refers to the existence of a function $K_r$ on $[-1,1]$ so that
$$K(x,y)=K_r(x \cdot y), \quad x,y \in S^1,$$
in which $\cdot$ is the usual inner product of $\mathbb{R}^2$.\ For a kernel $K$ on $S^1 \times S^1$, isotropy corresponds to the property
$$K((x,z),(y,w))=K_r(x \cdot y, z\cdot w), \quad x,y,z,w \in S^1,$$
in which the function $K_r$ has now domain $[-1,1]^2$.\ In both cases, we will call $K_r$ the {\em isotropic part} of the kernel $K$.\ Finally, the {\em positive definiteness} of a real kernel $K$ on an infinite set $X$ refers to the validity of the inequality
$$
\sum_{\mu,\nu=1}^n c_\mu c_\nu K(x_\mu, x_\nu) \geq 0,
$$
whenever $n$ is a positive integer, $x_1, x_2, \ldots, x_n$ are distinct points on $X$ and the $c_\mu$ are real scalars.\ The {\em strict positive definiteness} of $K$ demands both, its positive definiteness and that the inequalities above be strict whenever at least one of the $c_\mu$ is nonzero.\ We will apply these definitions to the cases in which either $X=S^1$ or $X=S^1 \times S^1$.

According to a result of Schoenberg in \cite{schoen}, a real, continuous and isotropic kernel $K$ on $S^1$ is positive definite if, and only if, the isotropic part $K_r$ of $K$ has the form
$$
K_r(t)=\sum_{k=0}^\infty a_k P_k(t), \quad t \in [-1,1],
$$
in which all the $a_k$ are nonnegative, $P_k$ is the Tchebyshev polynomial (of first kind) of degree $k$ (see \cite{szego}), and $\sum_{k=0}^\infty a_k P_k(1)<\infty$.\ In the nineties, many attempts were made in order to deduce a similar characterization for the strict case (\cite{ron,sun}) but that only appeared in \cite{mene} (see also \cite{barbosa}): a kernel having Schoenberg's representation described above is strictly positive definite on $S^1$ if, and only if, the set $\{k:a_{|k|}>0\}$ intersects every full arithmetic progression in $\mathbb{Z}$.\ Both results described above extends to the complex setting, that is, to the case in which $S^1$ is replaced with the unit circle in $\mathbb{C}$, the positive definite kernel is allowed to assume complex values and the scalars $c_\mu$ are now complex numbers.\ This extension is also discussed in \cite{mene}.

Moving to $S^1 \times S^1$, a theorem proved in \cite{jean} includes a characterization for the positive definiteness of a real, continuous and isotropic kernel $K$ on $S^1 \times S^1$ as those having an isotropic part in the form
\begin{equation}\label{repre}
K_r(t,s)=\sum_{k,l=0}^\infty a_{k,l} P_k(t)P_l(s), \quad t,s \in [-1,1],
\end{equation}
where all the coefficients $a_{k,l}$ are nonnegative and $\sum_{k,l=0}^\infty a_{k,l} P_k(1)P_l(1)<\infty$.\ The results in this paper will converge to the following characterization for strict positive definiteness on $S^1 \times S^1$.

\begin{thm} \label{mainintro} Let $K$ be a real, continuous, isotropic and positive definite kernel on $S^1 \times S^1$.\ It is strictly positive definite if, and only if,
the set $\{(k,l): a_{|k|,|l|} >0\}$ extracted from the series representation for the isotropic part $K_r$ of $K$ intersects all the translations of every lattice of $\mathbb{Z}^2$.
\end{thm}

The paper proceeds as follows.\ In Section 2, we present a technical result that leads to an alternative formulation for the concept of strict positive definiteness and use it to ratify the necessity part of Theorem \ref{mainintro}.\ In Section 3, after presenting a series of three technical results, we finally state and prove a result that not only includes Theorem \ref{mainintro} but also provides an alternative characterization for the strict positive definiteness on $S^1 \times S^1$.\ In Section 4, we provide a different path one could follow in order to simplify the proof of Theorem \ref{repre}.\ In Section 5, we indicate how the main results of the paper can be extended to the complex setting.

\section{Strict positive definiteness on $S^1 \times S^1$}

This section is divided in two parts: in the first one we explore a little bit deeper the concept of strict positive definiteness on $S^1 \times S^1$.\ The outcome Proposition \ref{prop-equiv cAc=0} implies an obvious equivalence for that concept.\ In the second part, we use this alternative description in order to obtain a proof for the necessity half of Theorem \ref{mainintro}.

For a function $K_r$ representable as in (\ref{repre}), we will write
$$J_K:=\{(k,l)\in\mathbb{Z}_+^2: a_{k,l}>0\}.$$
It is an easy matter to verify that the strict positive definiteness of the kernel $K$ with isotropic part given by (\ref{repre}) depends upon $J_K$ only and not on the actual values of the Fourier coefficients $a_{k,l}$.\
For distinct points $(x_1,w_1),(x_2,w_2), \ldots,(x_n,w_n)$ on $S^1 \times S^1$, we will write $A=(A_{\mu\nu})$, in which
$$
A_{\mu\nu}=K_r(x_\mu \cdot x_\nu, w_\mu \cdot w_\nu), \quad \mu,\nu=1,2,\ldots,n.
$$
We will also represent the points above in polar form:
$$
x_\mu=(\cos \t_\mu, \sin \t_\mu), \quad w_\mu=(\cos \p_\mu, \sin \p_\mu), \quad \t_\mu,\p_\mu \in [0,2\pi),\quad \mu=1,2,\ldots,n.
$$

\begin{prop} \label{prop-equiv cAc=0}
Let $K$ be a nonzero, real, continuous, isotropic and positive definite kernel on $S^1 \times S^1$.\ For distinct points $(x_1,w_1),(x_2,w_2), \ldots,(x_n,w_n)$ on $S^1 \times S^1$ and a column
vector $c=(c_\mu)$ in $\mathbb{R}^n$, the following statements are equivalent:\\
$(i)$ $c^{t}Ac=0$;\\
$(ii)$ The double equality
$$\sum_{\mu=1}^{n}c_{\mu}e^{ik\t_{\mu}}e^{il\p_{\mu}}=\sum_{\mu=1}^{n}c_{\mu}e^{ik\t_{\mu}}e^{-il\p_{\mu}}=0$$
holds for all $(k,l)\in J_K$.
\end{prop}
\pf The normalization one decides to adopt for the Tchebyshev polynomials is of no importance in this paper.\ So, we will write
$$
P_k(x_\mu \cdot x_\nu)=\frac{2}{k} \cos k(\t_\mu -\t_\nu), \quad k>0, \quad x_\mu \in S^1, \quad \mu=1,2,\ldots,n,
$$
while
$P_0(x_\mu \cdot x_\nu)=1$, $\mu=1,2,\ldots,n$.\ Now, computing each expression
$$
\sum_{\mu,\nu=1}^{n} c_\mu c_\nu P_k(x_\mu \cdot x_\nu)P_l(w_\mu \cdot w_\nu)
$$
in the cases when both $k$ and $l$ are nonzero and when at least one of them is zero, it is not hard to see that
the equality $c^{t}Ac=0$ corresponds to
$$
\left|\sum_{\mu=1}^{n}c_{\mu}e^{ik\t_{\mu}}e^{il\p_{\mu}}\right|^2+
\left|\sum_{\mu=1}^{n}c_{\mu}e^{ik\t_{\mu}}e^{-il\p_{\mu}}\right|^2=0, \quad (k,l)\in J_K.
$$
But that is equivalent to the statement in $(ii)$.\eop\vspace*{3mm}

From now on, we will deal with subgroups of $\mathbb{Z}^2$ and their translations.\ The subgroups of $\mathbb{Z}^2$ can be classified as follows.

\begin{lem}\label{class}
A nontrivial subgroup of $\mathbb{Z}^2$ belongs to one of the following categories:\\
$(i)$ $(0,b)\mathbb{Z}:=\{(0,pb): p \in \mathbb{Z}\}$, $b>0$; \\
$(ii)$ $(a,b)\mathbb{Z}:=\{(pa, pb): p \in \mathbb{Z}\}$, $a>0$;\\
$(iii)$ $(a,b)\mathbb{Z}+(0,d)\mathbb{Z}:=\{(pa,pb+qd): p,q \in \mathbb{Z}\}$, $a,d>0$.
\end{lem}

We advise the reader that there are different ways to describe the subgroups of $\mathbb{Z}^2$ (for instance, the one presented in \cite{pinkus} is slightly different and quite elegant).\ The subgroups that fit into Lemma \ref{class}-$(iii)$ will be called {\em lattices}.\ The set of lattices of $\mathbb{Z}^2$ encompasses all the subgroups of rank 2.\ If $ad=1$, then a lattice becomes the whole $\mathbb{Z}^2$, otherwise it is a proper subgroup of $\mathbb{Z}^2$.\ The lattices having the form
$$
(a\mathbb{Z},b\mathbb{Z}):=\{(pa, qb): q,p\in \mathbb{Z}\},\quad a,b>0,
$$
will be called {\em rectangular lattices} of $\mathbb{Z}^2$.\ By {\em translates of subgroups} of $\mathbb{Z}^2$, we will mean sets of the form
$(j,j')+S$, in which $(j,j')$ is a fixed element of $\mathbb{Z}^2$ and $S$ is a subgroup of $\mathbb{Z}^2$.

The next result is very close to the necessity part of Theorem \ref{mainintro}.

\begin{thm} \label{juru} Let $K$ be a real, continuous, isotropic and positive definite kernel on $S^1 \times S^1$.\
If $K$ is strictly positive definite, then the set $\{(k,l): (|k|,|l|)\in J_K\}$                                                                                                                                                                                                                                    intersects all the translations of each rectangular lattice of $\mathbb{Z}^2$.
\end{thm}
\pf Assume $K$ is strictly positive definite and write $S=(a\mathbb{Z}, b\mathbb{Z})$ with $a,b>0$.\ We will show that $$\{(k,l): (|k|,|l|)\in J_K\}$$ intersects $(j,j')+S$, whenever $j\in \{0,1,\ldots, a-1\}$ and $j'\in \{0,1,\ldots, b-1\}$.\ There is nothing to prove if $a=b=1$.\ In the other cases, we will assume that
$$\{(k,l): (|k|,|l|)\in J_K\} \cap (j+a\mathbb{Z}, j'+b\mathbb{Z})=\emptyset,$$
and will reach a contradiction.\ In the case in which $a=1$ and $b\geq 2$, the assumption on $J_K$ implies that $l-j',-l-j' \not \in b\mathbb{Z}$, whenever $(k,l) \in J_K$.\ In particular,
$$\sum_{\mu=1}^b\left(e^{i2\pi \mu /b}\right)^{l-j'}=\sum_{\mu=1}^b \left(e^{i2\pi \mu /b}\right)^{-l-j'}=0,$$
and, consequently,
$$\sum_{\mu=1}^b \left[\mbox{Re\,}\left(e^{i2\pi \mu j'/b}\right)\right] \left(e^{i2\pi \mu/b}\right)^l=0, \quad (k,l) \in J_K.$$
The real scalars $c_\mu:=\mbox{Re\,}(e^{i2\pi \mu j'/b})$, $\mu=1,2,\ldots,b$, are not all zero and the points
$$(x_\mu,w_\mu)=(1, e^{i2\pi \mu/b}), \quad \mu=1,2,\ldots,b,$$ are distinct in
$S^1\times S^1$.\ Thus, under the light of Proposition \ref{prop-equiv cAc=0}, we have a contradiction with the strict positive
definiteness of $K$.\ The case in which $a\geq 2$ and $b=1$ is similar.\ To conclude the proof, we now assume $a,b\geq 2$ and adapt the procedure employed in the first case.\ If
$(k,l) \not \in (j+a\mathbb{Z}, j'+b\mathbb{Z})$, then
either $k-j \not \in a\mathbb{Z}$ or $l-j' \not \in b\mathbb{Z}$.\ Hence, we may conclude that
$$\sum_{\mu=1}^a \left(e^{i2\pi \mu /a}\right)^{k-j} \sum_{\nu=1}^b\left(e^{i2\pi \nu /b}\right)^{l-j'}=0,$$
that is,
$$\sum_{\mu=1}^a \sum_{\nu=1}^b e^{-i2\pi \mu j /a}e^{-i2\pi \nu j'/b} \left(e^{i2\pi \mu /a}\right)^k \left(e^{i2\pi \nu /b}\right)^l=0.
$$
Repeating the argument with the assumption $(-k,-l) \not \in (j+a\mathbb{Z}, j'+b\mathbb{Z}),$
we conclude that
$$\sum_{\mu=1}^a \sum_{\nu=1}^b e^{i2\pi \mu j /a}e^{i2\pi \nu j'/b} \left(e^{i2\pi \mu /a}\right)^k \left(e^{i2\pi \nu /b}\right)^l=0.
$$
Thus, since $(k,l)$ is arbitrary,
$$\sum_{\mu=1}^a \sum_{\nu=1}^b \left[\mbox{Re\,}\left(e^{i2\pi \mu j /a}e^{i2\pi \nu j'/b}\right)\right] \left(e^{i2\pi \mu /a}\right)^k \left(e^{i2\pi \nu /b}\right)^l=0, \quad (k,l) \in J_K.
$$
By an analogous procedure, now taking into account that $(-k,l),(k,-l) \not \in (j+a\mathbb{Z}, j'+b\mathbb{Z})$, the conclusion is
$$\sum_{\mu=1}^a \sum_{\nu=1}^b \left[\mbox{Re\,}\left(e^{i2\pi \mu j /a}e^{i2\pi \nu j'/b}\right)\right] \left(e^{i2\pi \mu /a}\right)^k \left(e^{-i2\pi \nu /b}\right)^l=0, \quad (k,l) \in J_K.
$$
Therefore, since the numbers $\mbox{Re\,}\left[e^{i2\pi \mu j /a}e^{i2\pi \nu j'/b}\right]$ are not all zero and the $ab$ points
$$(x_\mu,w_\nu)=(e^{i2\pi \mu/a}, e^{i2\pi \nu/b}), \quad \mu=1,2,\ldots,a,\quad \nu=1,2,\ldots, b,$$ are distinct in $S^1 \times S^1$, we have reached a contradiction once again.\eop

In order to ratify the necessity part of Theorem \ref{mainintro}, the following lemma becomes handy.

\begin{lem} \label{tech} The lattice $L=(a,b)\mathbb{Z}+(0,d)\mathbb{Z}$, $a,d>0$, can be decomposed in the form
$$L=\bigcup_{(j,j')\in A}\left[(j,j')+(ad\mathbb{Z}, ad\mathbb{Z})\right],$$
in which $A=L\cap \{(\alpha,\beta) \in \mathbb{Z}^2: 0\leq \alpha,\beta < ad\}$.
\end{lem}
\pf For $p,q \in \mathbb{Z}$, we can certainly write
$$(pa,pb+qd)=(j+\alpha ad, j'+ \beta ad),$$
in which $j,j'\in \{0,1,\ldots,ad-1\}$.\ Since
$$(j,j')=(p-\alpha d)(a,b)+(\alpha b +q-\beta a)(0,d),$$
it is clear that $(j,j') \in L$.\ These arguments show that $L$ is a subset of the union quoted in the statement of the lemma.\ As for the reverse inclusion,
first observe that if $\alpha,\beta \in \mathbb{Z}$, we have that
$$(\alpha ad,\beta ad)=d\alpha (a, b) +(a\beta-b\alpha)(0,d) \in L.$$
Since $L$ is a subgroup of $\mathbb{Z}^2$, $(j,j')+(\alpha ad,\beta ad) \in L$ whenever $(j,j')\in A$.\eop\vspace*{3mm}

The theorem below is now evident.

\begin{thm} \label{main2}
Let $K$ be a real, continuous, isotropic and positive definite kernel on $S^1 \times S^1$.\ If $K$ is strictly positive definite, then the set $\{(k,l): (|k|,|l|)\in J_K\}$ intersects
all the translations of each lattice in $\mathbb{Z}^2$.
\end{thm}

\section{Sufficiency via M. Laurent's theorem}

In this section, we will prove that the necessary conditions for the strict positive definiteness of a continuous, isotropic and positive definite kernel on $S^1 \times S^1$ presented in Theorem \ref{juru} are also sufficient.

We begin recalling an elementary bi-dimensional version of the Skolem-Mahler-Lech Theorem due to M. Laurent (\cite{laurent,laurent1}).\ The original Skolem-Mahler-Lech Theorem is discussed in details in \cite{everest}.\ This very same bi-dimensional version was used in \cite{pinkus} in order to characterize certain strictly positive definite kernels on complex Hilbert spaces.

\begin{thm} \label{laurent}
Let $\{(x_1,w_1), (x_2,w_2), \ldots, (x_n,w_n)\}$ be a subset of $\mathbb{C}^2$.\ For $n$ complex numbers $c_1,c_2, \ldots, c_n$, define a double sequence $\{b_{k,l}:k,l \in \mathbb{Z}\}$
through the formula
$$b_{k,l}:=\sum_{\mu=1}^{n}c_\mu x_\mu^k w_\mu^l, \quad k,l\in \mathbb{Z}.$$
Then, the set $\{(k,l): b_{k,l}=0\}$ is the union of a finite number of translates of subgroups of $\mathbb{Z}^2$.
\end{thm}

The technical lemma below adds to Theorem \ref{laurent} when the points are distinct and belong to $\Omega_2 \times \Omega_2$, in which $\Omega_2$ is the unit circle in $\mathbb{C}$.

\begin{lem} \label{notallzero} Let $(x_1,w_1), (x_2,w_2), \ldots, (x_n,w_n)$ be distinct points in $\Omega_2 \times \Omega_2$.\ For complex number $c_1,c_2, \ldots, c_n$, define
$$b_{k,l}=\sum_{\mu=1}^{n}c_\mu x_\mu^k w_\mu^l, \quad k,l \in \mathbb{Z}.$$
If $\{(k,l): b_{k,l}=0\}=\mathbb{Z}^2$, then all the $c_\mu$ are zero.
\end{lem}
\pf We will write the components of the points in polar form $x_\mu=e^{i\t_\mu}$, $w_\mu=e^{i\p_\mu}$, $\mu=1,2,\ldots,n$, and will assume, as we can, that the $n$ points $(\t_1,\p_1), (\t_2, \p_2), \ldots, (\t_n, \p_n)$ are distinct in $[0,2\pi)^2$.\ Choose $\alpha,\beta \in \mathbb{Z}$ in such a way that all the elements in the set
$$
\left\{\alpha\frac{\t_\mu-\t_\nu}{2\pi}+\beta \frac{\p_\mu-\p_\nu}{2\pi}: \mu,\nu=1,2,\ldots,n; \mu \neq \nu\right\}
$$
are nonzero.\ Next, pick $\gamma\in \mathbb{Z}_+$ arbitrarily large so that
$$\left\{\frac{\alpha}{\gamma}\frac{\t_\mu-\t_\nu}{2\pi}+\frac{\beta}{\gamma} \frac{\p_\mu-\p_\nu}{2\pi}: \mu,\nu=1,2,\ldots,n; \mu \neq \nu\right\}\subset (-1,1)\setminus \{0\}.$$
For each pair $(\mu,\nu)$, $\mu \neq \nu$, for which
$$\frac{\alpha}{\gamma} \frac{\t_\mu-\t_\nu}{2\pi}+\frac{\beta}{\gamma} \frac{\p_\mu-\p_\nu}{2\pi} \in \mathbb{Q},$$
let $p_{\mu\nu}$ be a positive integer $>\gamma$ satisfying
$$p_{\mu\nu}\left(\frac{\alpha}{\gamma} \frac{\t_\mu-\t_\nu}{2\pi}+\frac{\beta}{\gamma} \frac{\p_\mu-\p_\nu}{2\pi} \right)\in \mathbb{Z}.$$
Finally, select an integer $q$ so that $q$ is greater then all the $p_{\mu\nu}$ and each set $\{q,p_{\mu\nu}\}$ is coprime.\ If $\{(k,l): b_{k,l}=0\}=\mathbb{Z}^2$, then we may infer that
$$
\sum_{\mu=1}^n c_\mu e^{i\t_\mu k}e^{i\p_\mu l}=0, \quad (k,l)=(0,0),(\alpha q,\beta q),(2\alpha q,2\beta q),\ldots, ((n-1)\alpha q,(n-1)\beta q).
$$
The matrix of the system above has $\mu\nu$-entries given by
$$
\left[e^{i(\alpha \t_\mu +\beta \p_\mu)q}\right]^{\nu},\quad \mu,\nu=1,2,\ldots,n,
$$
and, consequently, it is a Vandermonde matrix.\ So, the proof of the lemma will be complete as long as we show that the $n$ points $e^{i(\alpha \t_\mu +\beta \p_\mu)q}$, $\mu=1,2,\ldots, n$,
are distinct.\ But, for $\mu \neq \nu$,
$$e^{i(\alpha \t_\mu +\beta \p_\mu)q}= e^{i(\alpha \t_\nu +\beta \p_\nu)q}$$
if, and only if,
$$
q \left(\alpha\frac{\t_\mu -\t_\nu}{2\pi} +\beta\frac{\p_\mu -\p_\nu}{2\pi}\right)\in \mathbb{Z}.
$$
If all the numbers
$$\left(\alpha\frac{\t_\mu -\t_\nu}{2\pi} +\beta\frac{\p_\mu -\p_\nu}{2\pi}\right),\quad \mu \neq \nu,$$
are irrational, we are done.\ Otherwise, there would be integers $j$ and $j'$ such that
$$\frac{j}{\gamma q}=\frac{ j'}{p_{\mu\nu}}\in (-1,1)\setminus \{0\}.$$
for some pair $(\mu,\nu)$, $\mu \neq \nu$.\ Since $\{q,p_{\mu\nu}\}$ is coprime, then $p_{\mu\nu}$ would divide $\gamma$, contradicting our choice of $p_{\mu\nu}$.\eop\vspace*{3mm}

The next result reveals that if a proper subset $A$ of $\mathbb{Z}^2$ is a finite union of translates of subgroups of $\mathbb{Z}^2$, then there exists a rectangular lattice $H$ of $\mathbb{Z}^2$ and $(j,j') \in \mathbb{Z}^2$ so that $[(j,j')+H]\cap A=\emptyset$.

\begin{lem} \label{elattice} Let $A$ be a proper subset of $\mathbb{Z}^2$.\ If $A$ is a finite union of translates of subgroups of $\mathbb{Z}^2$ and $(j,j')\in \mathbb{Z}^2\setminus A$, then there exists a rectangular lattice $H$ of $\mathbb{Z}^2$ such that $(j,j')+H \subset \mathbb{Z}^2\setminus A$.
\end{lem}
\pf If $A$ is a finite union of translates of subgroups of $\mathbb{Z}^2$, we can write
$$A=F \cup [(j_1,j_1') +G_1] \cup [(j_2,j_2') +G_2] \cup \cdots \cup [(j_r,j_r') +G_r]$$
in which $F$ is a finite (possibly empty) subset of $\mathbb{Z}^2$, $(j_1,j_1'), (j_2,j_2'), \ldots, (j_r,j_r')\in \mathbb{Z}^2$ and $G_1, G_2, \ldots, G_r$ are nontrivial subgroups of $\mathbb{Z}^2$.\ It suffices to prove the lemma in the case in which $F=\emptyset$.\ Indeed, if a solution $(j,j')+H$ is available for that case, we can pick a convenient subgroup $H_1$ of $H$ so that $(j,j')+H_1$ avoids all the elements of $F$.\ So, assume that $F=\emptyset$ and fix $(j,j') \in \mathbb{Z}^2 \setminus A$.\ We can assume all the $G_i$ have rank 2.\ Indeed, if $G_i$ has rank 1 for some $i$, we can pick $(\alpha,\beta) \in \mathbb{Z}^2$ such that
$$(\alpha , \beta)\mathbb{Z} \cap [(j_i-j,j_i'-j')+G_i]=\emptyset.$$
Hence,
$$[(j,j')+(\alpha , \beta)\mathbb{Z}] \cap [(j_i,j_i')+G_i]=\emptyset,$$
and, therefore,
$$(j,j') \not \in (j_i,j_i')+(\alpha , \beta)\mathbb{Z} +G_i.$$
In particular, $(\alpha , \beta)\mathbb{Z} +G_i$ is a subgroup of rank 2 and we can replace $(j_i,j_i)+G_i$ with
$(j_i,j_i')+(\alpha , \beta)\mathbb{Z} +G_i$ in the union decomposition for $A$ keeping $(j,j')$ in $\mathbb{Z}^2\setminus A$.\ If all the $G_i$ have rank 2, the proof of the lemma proceeds as follows.\ Let $m_i$ be the index of $G_i$ in $\mathbb{Z}^2$, $i=1,2,\ldots,r$, and pick a common multiple $m$ of all the $m_i$.\ The subgroup $(m\mathbb{Z},m\mathbb{Z})$ is a rectangular lattice and, by the definition of index of a subgroup, it follows that
$$(m\mathbb{Z},m\mathbb{Z}) \subset G_i \quad i=1,2,\ldots,r.$$
In particular, $$[(j,j')+(m\mathbb{Z},m\mathbb{Z})]\cap G_i=\emptyset, \quad i=1,2,\ldots,r,$$
and, consequently, $(j,j')+(m\mathbb{Z},m\mathbb{Z}) \subset \mathbb{Z}^2\setminus A$. \eop \vspace*{3mm}

Next, we prove a refinement of the sufficiency part of Theorem \ref{mainintro}.

\begin{thm} \label{sufic} Let $K$ be a real, continuous, isotropic and positive definite kernel on $S^1 \times S^1$.\ If $\{(k,l): (|k|,|l|)\in J_K\}$ intersects all the translations of each rectangular lattice of $\mathbb{Z}^2$, then $K$ is strictly positive definite.
\end{thm}
\pf Let us assume that $\{(k,l): (|k|,|l|)\in J_K\}$ intersects all the translations of each rectangular lattice of $\mathbb{Z}^2$.\ For a fixed $n\geq 2$, $n$ distinct points $(\t_1,\p_1), (\t_2, \p_2), \ldots, (\t_n, \p_n)$ in $[0,2\pi)^2$ and $c_1, c_2, \ldots, c_n$ real numbers, not all zero, we intend to show that either $\sum_{\mu=1}^n c_\mu e^{i\t_\mu k}e^{-i\p_\mu l} \neq 0$ or $\sum_{\mu=1}^n c_\mu e^{i\t_\mu k}e^{i\p_\mu l}\neq 0$, for some $(k,l) \in J_K$.\ A help of Proposition \ref{prop-equiv cAc=0} will lead to the strict positive definiteness of $K$.\ In order to achieve the conclusion mentioned above, define
$$b_{k,l}=\sum_{\mu=1}^n c_\mu e^{i\t_\mu k}e^{i\p_\mu l}, \quad k,l\in \mathbb{Z}.$$
On one hand, Lemma \ref{notallzero} and the fact that at least one $c_\mu$ is nonzero imply that $\{(k,l): b_{k,l}=0\}\neq \mathbb{Z}^2$.\
Theorem \ref{laurent} asserts that $\{(k,l): b_{k,l}=0\}$ is the union of a finite number of translations of subgroups of $\mathbb{Z}^2$ while Lemma \ref{elattice} guarantees the existence of a rectangular lattice of $\mathbb{Z}^2$, a translation of which belongs to $\mathbb{Z}^2 \setminus \{(k,l): b_{k,l}=0\}$.\ Thus, due to our assumption on $\{(k,l) : (|k|,|l|) \in J_K\}$, we immediately have that
$$
\{(k,l) : (|k|,|l|) \in J_K\}\not \subset \{(k,l): b_{k,l}=0\}.
$$
Therefore, there must exist at least one pair $(k,l)$ in $\{(k,l): (|k|,|l|)\in J_K\}$ for which
$$\sum_{\mu=1}^n c_\mu e^{i\t_\mu k}e^{i\p_\mu l}\neq 0.$$
Since the $c_\mu$ are real, the result follows.
\eop\vspace*{3mm}

The following characterizations are now evident.

\begin{thm}
Let $K$ be a real, continuous, isotropic and positive definite kernel on $S^1 \times S^1$.\ The following assertions are equivalent:\\
$(i)$ $K$ is strictly positive definite;\\
$(ii)$ $\{(k,l): (|k|,|l|)\in J_K\}$ intersects all the translations of each lattice in $\mathbb{Z}^2$;\\
$(iii)$ $\{(k,l): (|k|,|l|)\in J_K\}$ intersects all the translations of each rectangular lattice  in $\mathbb{Z}^2$.
\end{thm}

To conclude the section, we present an alternative guise for the previous characterizations.\ That will require the following technical, but elementary lemma.

\begin{lem}
Let $K$ be a subset of $\mathbb{Z}^2$ that intersects all the translations of each lattice in $\mathbb{Z}^2$.\ Then, every such intersection is an infinite set.
\end{lem}
\pf Let $L=(j,j')+(a,b)\mathbb{Z}+(0,d)\mathbb{Z}$, $a,d>0$, and assume that $K\cap L$ is finite.\ Write
$$(j+p_1a,j'+p_1b+q_1d), (j+p_2a,j'+p_2b+q_2d), \ldots, (j+p_na,j'+p_nb+q_nd),$$
to denote the elements in the intersection and define
$$
p:=\max\{|p_1|,|p_2|,\ldots,|p_n|\} \quad \mbox{and} \quad q:=\max\{|q_1|,|q_2|,\ldots,|q_n|\}.
$$
We will reach a contradiction, analyzing four different cases.\\
{\em Case1.} $p=q=0$: The intersection contains just one element, $(j,j')$.\ We now look at the translation
$$
L':=(j+2a,j'+2b)+(3a,3b)\mathbb{Z} +(0,d)\mathbb{Z}\subset L
$$
of the sublattice $(3a,3b)\mathbb{Z} +(0,d)\mathbb{Z}$ of $(a,b)\mathbb{Z}+(0,d)\mathbb{Z}$.\ If $(j,j')\in L'$, then
$$
\left\{
\begin{array}{c}
  j+2a+3ar = j \\
  j'+2b+3br+ds = j'
\end{array}
\right.
$$
for some $r,s\in\mathbb{Z}$.\ But, since $a(3r+2)\neq0$, $r\in\mathbb{Z}$, this is impossible.\ In particular, $K\cap L'=\emptyset$, a contradiction to our basic assumption.\\
{\em Case 2.} {\em $p=0$ and $q>0$}: Here we consider the sublattice $(a,b)\mathbb{Z} +(0,{2(2q+1)}d)\mathbb{Z}$ of $(a,b)\mathbb{Z}+(0,d)\mathbb{Z}$ and look at its translation
$$
L'':=(j,j'+2qd) + (a,b)\mathbb{Z} +(0,2(2q+1)d)\mathbb{Z} \subset L.
$$
If $(j+ra,j'+2qd+rb+2s(2q+1)d)=(j,j'+q_\mu d)$ for some $\mu \in \{1,2,\ldots,n\}$ and $r,s \in \mathbb{Z}$, then
$$
\left\{
\begin{array}{c}
  ra=0\\
  2qd+rb+2s(2q+1)d=q_\mu d
\end{array}
\right.
$$
and, consequently, $2q+2s(2q+1)=q_\mu$.\ However, due to the definition of $q$, no integer $s$ can satisfy the previous equality.\ Thus, $L''\cap K=\emptyset$, another contradiction.\\
{\em Case 3.} {\em $p>0$ and $q=0$}: Its is similar to the previous case.\\
{\em Case 4.} $p,q>0$: Here we consider the sublattice $(2(2p+1)a,2(2p+1)b)\mathbb{Z} +(0,qd)\mathbb{Z}$ of $(a,b)\mathbb{Z}+(0,d)\mathbb{Z}$ and its translation
$$L''':= (j+2pa,j'+2pb)+(2(2p+1)a,2(2p+1)b)\mathbb{Z} +(0,qd)\mathbb{Z} \subset L.$$
If
$$
(j+2pa+2r(2p+1)a, j'+2pb+2r(2p+1)b + sqd)= (j+p_\mu a,j'+p_\mu b+q_\mu d)
$$
for some $\mu \in \{1,2,\ldots,n\}$ and $r,s \in \mathbb{Z}$, we will have that
$2p+2r(2p+1)=p_\mu$.\ As in Case 2, we can deduce that $L'''\cap K=\emptyset$, a contradiction to our initial assumption on $K$.\eop

It is now clear that the following additional result holds.

\begin{thm}
Let $K$ be a real, continuous, isotropic and positive definite kernel on $S^1 \times S^1$.\ The following assertions are equivalent:\\
$(i)$ $K$ is strictly positive definite;\\
$(ii)$ $\{(k,l): (|k|,|l|)\in J_K\}$ intersects all the translations of each lattice in $\mathbb{Z}^2$ infinitely many times;\\
$(iii)$ $\{(k,l): (|k|,|l|)\in J_K\}$ intersects all the translations of each rectangular lattice  in $\mathbb{Z}^2$ infinitely many times.
\end{thm}

\section{Appendix: a direct proof for Theorem \ref{main2}}

Here, we include a self-contained direct proof for  Theorem \ref{main2}.\ A few passages incorporates new arguments.

\pf Let $L=(j,j')+(a,b)\mathbb{Z}+(0,d)\mathbb{Z}$, $a,d\geq1$, be a lattice of $\mathbb{Z}^2$.\ If $d=1$, then we can assume that $a>1$ because, otherwise, $L=\mathbb{Z}^2$ and the the conclusion of the Theorem holds trivially.\ We will suppose that $K$ is strictly positive definite and that
$\{(k,l): (|k|,|l|)\in J_K\} \cap L=\emptyset$ and will reach a contradiction.\ For each $(k,l)\in J_K$, we consider all possible cases attached to our
assumption.\ Since $(k,l) \not \in L$, either $k-j \not \in a\mathbb{Z}$ or
$$
k-j \in a\mathbb{Z} \quad \mbox{and } \quad l -j' -(k-j)b/a  \not \in d\mathbb{Z}.
$$
In the first possibility, we have that $a\geq 2$ and we can conclude that
$$
\sum_{\mu=1}^a e^{-i2\pi \mu j/a} e^{i 2 \pi \mu k/a}=\sum_{\mu=1}^{a} (e^{i2\pi \mu/a})^{k-j}=0.
$$
In the other, $d\geq 2$ and the conclusion is
$$
\sum_{\nu=1}^d  \left[e^{i2\pi \nu /d} \right]^{l-j'-(k-j)b/a}=0.
$$
The final outcome is
$$\sum_{\mu=1}^a  \left[e^{i2\pi \mu /a} \right]^{k-j} \sum_{\nu=1}^d  \left[e^{i2\pi \nu /d} \right]^{l-j'-(k-j)b/a}=0, \quad (k,l) \in J_K.$$
A similar argument now using the fact that $(-k,-l) \not \in L$, leads to
$$\sum_{\mu=1}^a  \left[e^{i2\pi \mu /a} \right]^{-k-j} \sum_{\nu=1}^d  \left[e^{i2\pi \nu /d} \right]^{-l-j'+(k+j)b/a}=0, \quad (k,l) \in J_K.$$
These two pieces of information can be put together to ensure the equality
$$
\sum_{\mu=1}^a \sum_{\nu=1}^d \left[\mbox{Re\,}( e^{-i2\pi \mu j/a} e^{-i2\pi \nu j'/d} e^{i 2 \pi \nu j b/ad})\right] e^{i2\pi k(\mu-\nu b/d)/a} e^{i2\pi l \nu /d}= 0, \quad (k,l) \in J_K,$$
with $a\geq 2$ or $d\geq 2$.\ Repeating the procedure above in the cases $(-k,l) \not \in L$ and $(k,-l) \not \in L$, we can conclude that
$$
\sum_{\mu=1}^a \sum_{\nu=1}^d \left[\mbox{Re\,}( e^{-i2\pi \mu j/a} e^{-i2\pi \nu j'/d} e^{i 2 \pi \nu j b/ad})\right] e^{i2\pi k(\mu-\nu b/d)/a} e^{-i2\pi l \nu /d}=
 0, \quad (k,l) \in J_K,$$
under the same setting for $a$ and $d$.
The angles
$$
\theta_{\mu\nu}= \dfrac{2\pi\left(\mu-\nu b/d\right)}{a}, \quad \p_\nu=\dfrac{2\pi  \nu}{d},\quad \mu=1,2,\ldots,a, \quad \nu=1,2,\ldots,d,
$$
define $ad$ distinct points
$$
x_{\mu\nu}=(\cos\theta_{\mu\nu},\sin\theta_{\mu\nu}), \quad w_\nu=(\cos\p_\nu,\sin\phi_\nu),\quad \mu=1,2,\ldots,a, \quad \nu=1,2,\ldots,d.
$$
in $S^1 \times S^1$.\ Indeed, if $x_{\mu_1\nu}=x_{\mu_2\nu}$, for a fixed $\nu \in \{1,2,\ldots, d\}$ and $\mu_1, \mu_2 \in \{1,2,\ldots,a\}$, it is promptly seen that
$\t_{\mu_1\nu}=\t_{\mu_2\nu}$, that is, $\mu_1=\mu_2$.\ On the other hand, at least one of the numbers
$$r_{\mu\nu}:=\mbox{Re\,}( e^{-i2\pi \mu j/a} e^{-i2\pi \nu j'/d} e^{i 2 \pi \nu j b/ad}), \quad \mu=1,2,\ldots, a,\quad \nu=1,2,\ldots, d,$$
is not zero.\ Indeed, this is obvious if $j=0$.\ If $j\neq 0$ and $d>1$, then one of the three numbers $r_{11}$, $r_{2,1}$ and $r_{1d}$ must be nonzero.\ Otherwise, all three numbers
$$4\left(\frac{j}{a}+\frac{j'}{d}-\frac{j b}{ad}\right),\quad 4\left(\frac{2j}{a}+\frac{j'}{d}-\frac{jb}{ad}\right), \quad \mbox{and} \quad 4\left(\frac{j}{a}+j'-\frac{j b}{a}\right)$$
belong to $1+2\mathbb{Z}$.\ In particular,
$$4 j \in 2a\mathbb{Z},\quad \mbox{and}\quad 4(j+j'a-jb) \in a(1+2\mathbb{Z})$$
what generates a contradiction (even=odd).\
If $j\neq 0$ and $d=1$, then $a\geq 2$ and there is just one value for $\nu$, namely $\nu=1$.\ In this case, the numbers $r_{\mu\nu}$ become
$$r_{\mu}:= \mbox{Re\,}[ e^{i 2\pi j(b-\mu)/a}],\quad \mu=1,2,\ldots, a.$$
If $a=2$, it is easily seen that  $r_1\neq 0$.\ If $a>2$, a calculation similar to another one made above yields that one of the three numbers $r_1$, $r_2$ and $r_a$ is nonzero.
Recalling Proposition \ref{prop-equiv cAc=0} once again, we now have reached a contradiction with the strict positive definiteness
of $K$.\eop

\section{Another appendix: the complex circle}

All the major results demonstrated in this paper can be adapted to hold for positive definiteness on the complex circle $\Omega_2$.\ In that case, we replace $S^1$ with $\Omega_2$, we allow the kernels to assume complex values and the scalars $c_\mu$ in the definition of positive definiteness can be complex numbers (the quadratic form in the definition of positive definiteness is Hermitian).\ We will sketch what these results are
and refer the interested reader to \cite{jean,mene} where the necessary adaptations for the proofs can be prospected from.

Let $K: \Omega_2 \times \Omega_2 \to \mathbb{C}$ be a continuous kernel and assume that
$$K((x,z),(y,w))=K_r(x\cdot y, z\cdot w),\quad x,y,z,w\in \Omega_2,$$
for some function $K_r : \Omega_2 \times \Omega_2 \to \mathbb{C}$, in which $\cdot$ is now the usual inner product of $\mathbb{C}$.\
It is positive definite if, and only if,
the function $K_r$ is of the form
$$
K_r(z,w)=\sum_{k,l\in\mathbb{Z}}a_{k,l}z^k w^l, \quad z,w \in \Omega_2,
$$
in which $a_{k,l} \geq 0$, $k,l \in \mathbb{Z}$ and $\sum_{k,l\in\mathbb{Z}}a_{k,l}<\infty$.\

Taking the above representation for granted, we can define $I_K:=\{(k,l): a_{k,l}>0\}$.\ For distinct points $(x_1,w_1),(x_2,w_2), \ldots,(x_n,w_n)$ on $\Omega_2 \times \Omega_2$ and a column vector $c$ in $\mathbb{C}^n$, the quadratic form $\overline{c}^{t}Ac=0$ corresponds to
$$\sum_{\mu=1}^{n}c_{\mu}e^{ik\t_{\mu}}e^{il\p_{\mu}}=0, \quad (k,l)\in I_K,$$
in which $\t_\mu$ and $\p_\mu$ are the arguments in the polar representation of $x_\mu$ and $w_\mu$ respectively.\ In particular, this reveals that the proofs we have developed in Sections 2 and 3 simplify in the present complex setting.

We close the paper stating what the main characterization for strict positive definiteness on $\Omega_2 \times \Omega_2$ becomes.

\begin{thm} Let $K$ be a continuous and positive definite kernel on $\Omega_2 \times \Omega_2$ as described above.\ It is strictly positive definite if, and only if,
 $I_K$ intersects all the translations of each rectangular lattice of $\mathbb{Z}^2$.
\end{thm}

\noindent {\bf Acknowledgement.} The arguments presented in the proof of Lemma \ref{elattice} are originally due to H. Borges and E. Tengan.\ We thank them for providing such elegant algebraic details.

\vspace*{15mm}

\noindent J. C. Guella, V. A. Menegatto and A. P. Peron \\
Departamento de
Matem\'atica,\\ ICMC-USP - S\~ao Carlos, Caixa Postal 668,\\
13560-970 S\~ao Carlos SP, Brasil\\ e-mails: jeanguella@gmail.com; menegatt@icmc.usp.br; apperon@icmc.usp.br

\end{document}